\documentclass{agtart_a}
\pdfoutput=1

\usepackage{enumerate, pb-diagram, graphicx}%

\title{Twisted Alexander polynomials of periodic knots}

\author{Jonathan A Hillman }
\givenname{Jonathan A}
\surname{Hillman}
\address{School of Mathematics and Statistics F07\\
University of Sydney\\\newline
NSW 2006\\
Australia}
\email{jonh@maths.usyd.edu.au}

\author{Charles Livingston}
\givenname{Charles}
\surname{Livingston}
\address{Department of Mathematics\\
Indiana University\\\newline
Bloomington IN 47405\\USA}
\email{livingst@indiana.edu}

\author{Swatee Naik}
\givenname{Swatee}
\surname{Naik}
\address{Department of Mathematics and Statistics\\
University of Nevada\\\newline
Reno NV 89557\\USA}
\email{naik@unr.edu}



\received{17 June  2005}
\revised{24 January 2006}
\accepted{26 January 2006}
\proposed{}
\seconded{}
\publishedonline{24 February 2006}
\published{24 February 2006}

\volumenumber{6}
\issuenumber{}
\publicationyear{2006}
\papernumber{5}
\startpage{145}
\endpage{169}

\doi{}
\MR{}
\Zbl{}

\AtBeginDocument{\let\bar\wbar}

\def\zz{\mathbb{Z}}
\def\qq{\mathbb{Q}}

\newcommand{\fig}[2] {\includegraphics[scale=#1]{\figdir/#2}}

\makeatletter
\def\cnewtheorem#1[#2]#3{\newtheorem{#1}{#3}
\expandafter\let\csname c@#1\endcsname\c@theorem}

\newtheorem{theorem}{Theorem}
\cnewtheorem{lemma}[theorem]{Lemma}
\cnewtheorem{corollary}[theorem]{Corollary}
\theoremstyle{definition}
\cnewtheorem{definition}[theorem]{Definition}

\makeatother  

\def\calc{\mathcal{C}}

\numberwithin{equation}{section}

\begin{document}

\begin{abstract} 
Murasugi discovered two criteria that must be satisfied by the
Alexander polynomial of a periodic knot.  We generalize these to the
case of twisted Alexander polynomials.  Examples demonstrate the
application of these new criteria, including to knots with trivial
Alexander polynomial, such as the two polynomial 1 knots with 11
crossings.

Hartley found a restrictive condition satisfied by the Alexander
polynomial of any freely periodic knot.  We generalize this result to
the twisted Alexander polynomial and illustrate the applicability of
this extension in cases in which Hartley's criterion does not apply.
\end{abstract}

\maketitle

\section{Introduction}  \label{sec:1}
The twisted Alexander polynomial of a knot, discovered independently
by Jiang--Wang \cite{j}, Lin \cite{li1} and Wada \cite{wa}, has seen a
growing number of applications to knot theory over the last ten years.
These range from the study of reversibility, mutation and concordance
of knots, Kirk--Livingston \cite{kl1,kl2} and Tamulis \cite{ta}, to
identifying fibered knots, Cha \cite{ch}, Goda--Morifuji \cite{gm} and
Goda--Kitano--Morifuji \cite{gkm}.  Most recently, Friedl and
Kim~\cite{fk} have demonstrated that the twisted polynomial is
sufficient to determine the genus and fibering properties of all prime
knots of 12 crossings and less.  Other literature on this invariant
includes Cogolludo--Florens \cite{cf}, Kitano \cite{ki},
Kitano--Suzuki \cite{ks}, Li--Xu \cite{ll} and Morifuji \cite{mo3,mo1}.
In this paper we extend the application of the twisted polynomial to
the study of knot periodicity.

A knot $K \subset S^3$ is called periodic, of period $q$, if there is
an order $q$ transformation of $S^3$ that leaves $K$ invariant and has
as fixed point set a circle $A$ disjoint from $K$.  The first major
result in the study of knot periodicity was the development of
algebraic obstructions by Trotter~\cite{tr}.  The central result
concerning algebraic properties of periodic knots is Murasugi's
theorem~\cite{mu} (see also Hillman \cite{hi1} and Sakuma \cite{sa}),
describing the Alexander polynomial of a periodic knot, $\Delta_K(t)$,
in terms of $\Delta_{\bar{L}}(t,s)$, the Alexander polynomial of the
quotient link, $\bar{L} = \bar{K} \cup \bar{A}$.  From this, two
conditions that must be satisfied by the Alexander polynomial of a
periodic knot follow.

\begin{theorem} \label{murthm1}
If $K$ is periodic of period $q$ with quotient $\bar{K}$, then there
is a polynomial $F(t,s) \in \zz[t^{\pm1}, s^{\pm 1}]$ such that $$
\Delta_K(t) = \Delta_{\bar{K}}(t) \prod_{i=1}^{q-1}F( t, \zeta_q^i),$$
where $\zeta_q$ is a primitive $q$--root of unity.
\end{theorem}

\noindent In a more precise statement of this theorem, it is seen that 
$F(t,s)  =  \Delta_{\bar{L}}(t,s)$,     and 
$F(t,1) =  \delta_\lambda (t)  \Delta_{\bar{K}}(t)$  where 
$\delta_\lambda(t) =  (1-t^\lambda)/(1-t)$ and  
$\lambda = \mbox{lk}(K, A)$.

\begin{theorem} \label{murthm2}  If $K$ is of  prime power period $q = p^r$ then 
$$\Delta_{K}(t) \equiv \Delta_{\bar{K}}(t)^{q} \delta_\lambda(t)^{q-1} 
\hskip.1in \mathrm{mod}   p.$$  
\end{theorem}

Since the publication of Murasugi \cite{mu}, our understanding of
periodicity of knots has greatly expanded.  In addition to a deeper
understanding of Murasugi's result (Burde \cite{bu}, Burde--Zieschang
\cite{bz}, Davis--Livingston \cite{dl1,dl2}, Hillman \cite{hi1,hi2}
and Sakuma \cite{sa}) new techniques have been applied, including
general knot polynomials, minimal surfaces and hyperbolic geometry;
examples of such work include Adams--Hildebrand--Weeks \cite{ahw},
Edmonds \cite{ed}, Naik \cite{na,na2}, Przytycki \cite{prz} Traczyk
\cite{tra} and Yokota \cite{yo}.  However, the Murasugi conditions
continue to be of central interest to the subject; in addition to
being easily applied to examples, their homological nature enables
their extension and application in other settings, including knots in
more general manifolds and higher dimensional knots.

In classical knot theory, in the case that abelian invariants do not 
resolve questions of interest, the next step has been  to move to 
nonabelian, for example dihedral, invariants of a knot.  Such 
invariants are placed in the context of homological invariants and 
combined with invariants of the infinite cyclic cover of a knot via 
the twisted Alexander polynomial.

Later we will give a detailed description of the twisted Alexander 
polynomial.  Briefly, given a link $L$  of $k$ components  and a 
representation  $\rho \co \pi_1(S^3 - L) \to GL_n(R)$ for some 
Noetherian unique factorization domain (such as a PID, ie,~principal ideal domain, or field)
$R$, there is a polynomial invariant $$\Delta(L, 
\rho)(t_1 , \ldots, t_k) \in R[t_1^{\pm 1}, \ldots, t_k^{\pm 1}],$$ well 
defined up to multiplication by a unit.  
We prove the following analogs of Murasugi's Theorems 1 and 2.

\begin{theorem}\label{tthm1} Suppose that $K$ is of period $q$, and $A$, $\bar{K}$, 
and $\bar{A}$ are as above.
Further, assume that $\bar{\rho} \co \pi_1(S^3 - \bar{K}) \to  GL_n(R)$, with $R = \zz$ or $\qq$, 
  and  $\rho$ is the lift of that
representation to $S^3 - K$.  Then there is a polynomial $F(t,s) \in R[t^{\pm 1},s^{\pm 1}]$ such that  $$\Delta_{K\!,  \rho}(t)  = 
\Delta_{\bar{K}\!, \bar{\rho}}(t) \prod_{i=1}^{q-1} F(t, \zeta_q^i).$$
\end{theorem}

\begin{theorem}\label{tthm2}  For a knot $K$ of period $q = p^r$, $p$ prime, and a  representation $\bar{\rho} \co \pi_1(S^3 - \bar{K}) \to GL_n(\zz / p\zz)$ with lift $\rho$,   if  $\Delta_{K\!,\rho}(t) \ne 0$, then 
 $$\Delta_{K\!,\rho}(t) = \Delta_{\bar{K}\!, \bar{\rho}}(t)^q 
(  \delta_{L ,  {\rho}}(t) /\Delta^0_{K\!,\rho}(t)) ^{q-1}.$$
\end{theorem}

\noindent Here $L = (K, A)$ and $ \delta_{L,  {\rho}}(t) \in ( \zz / p\zz)[t^{\pm 1}]$ 
is more 
complicated than in the classical case, but will be seen   in many 
cases to be  determined by    data about $K$ and $\rho$ alone. The 
term $\Delta^0$ in the denominator is related to twisted homology in 
dimension 0  and is easily computed in examples.  In the classical 
setting $\Delta^0$ is simply $(1-t)$ and $\delta_{L,\rho}(t) = 1 - t^\lambda$, 
so the classical $\delta_\lambda$ becomes our $\delta_{L ,  {\rho}}(t) /
\Delta^0_{K\!,\rho}(t)$.  Also, $F$ will be seen to be related to a 
twisted polynomial of $\bar{L}$.

The above theorems provide obstructions to periodicity 
when a representation of the knot group is necessarily a lift of that 
for a hypothesized quotient knot. 
When this fails to be the case, 
the twisted polynomials continue to
provide obstructions: 
there is a periodic action on the set of representations 
and different translates of a representation under this action must
yield the same twisted polynomial.

The study of freely periodic knots, those that are invariant under a 
free periodic transformation of $S^3$, has often been treated independently of the study of periodic knots.  We will see in the final section of this paper that the same techniques we use to study periodic knots apply to give results concerning the twisted Alexander polynomial of a freely periodic knot.  The main result concerning Alexander polynomials of freely periodic knots was proved by Hartley in~\cite{ha}.  We generalize his result to the twisted case as follows, where    $\bar{K}$ and $\bar{\rho}$ are defined similarly as above, though some care must be taken since the quotient  space is no longer $S^3$.  

\begin{theorem}\label{freethm} If $K$ is a freely periodic knot of period $q$ with quotient knot $\bar{K}$ having representation $\bar{\rho}$ that lifts to a representation $\rho$,  then 
$$\Delta_{K\!,\rho}(t^q) =   \prod_{i=0}^{q-1} \Delta_{\bar{K}\!,\bar{\rho}}(\zeta_q^i t).$$
\end{theorem}

{\bf Outline of Paper}\qua In the next section we will review the
general theory of twisted homology and its use in defining the twisted
Alexander modules of a knot.  \fullref{sec:3} considers torsion
invariants of modules, gives the definition of twisted Alexander
polynomials of knots and links, and describes techniques for their
computation.  In \fullref{sec:4} we relate the twisted Alexander
polynomial of a knot $K$ to that of a related link, $L = (K,A)$.
\fullref{sec:5} applies Shapiro's Lemma, which relates the twisted
homology of a space to that of a covering space, to give splitting
theorems for the twisted homology of links associated to periodic
knots.  The results of \fullref{sec:5} are applied in
\fullref{sectionproofrat} to prove \fullref{tthm1} above, the
generalization to twisted Alexander polynomials of Murasugi's result
on integral polynomials, which we stated in \fullref{murthm1}.
\fullref{sectiontwistp} proves \fullref{tthm2}, the twisted analog of
Murasugi's $\zz / p\zz$ result stated in \fullref{murthm2}.  In
\fullref{10_162} we provide basic examples illustrating these results
and in \fullref{poly1} we give more delicate examples using
knots with Alexander polynomial 1.  \fullref{sectionpres} discusses methods to
use the twisted polynomial to obstruct knot periodicity when
Theorems~\ref{tthm1} and~\ref{tthm2} do not apply because there are no
representations of the quotient knot that can be lifted.  Note that,
as all our examples are of hyperbolic knots, SnapPea, a program
written by Weeks~\cite{we} that computes the isometry group for a
hyperbolic knot complement provides periods of these knots. On the
other hand, the twisted polynomial approach is purely algebraic and
does not require a hyperbolic structure on the knot complement.  The
final section, \fullref{sec:11}, demonstrates that the techniques of this
paper extend to the study of freely periodic knots, giving an
extension of a classical result of Hartley~\cite{ha} to the case of
twisted Alexander polynomials.

\medskip{\bf Acknowledgments}\qua The authors wish to thank Jim Davis
and Paul Kirk for their helpful advice.  Computations of twisted
polynomials for many of the examples were run on the Beowulf cluster
at University of Nevada, Reno, maintained by Eric Olson; we thank Eric
for his help with the computations.  The first author received support
from Grey College and the University of Durham, through the Grey
College Mathematics Fellowship.  The second author received funding
from the {NSF}.  The third author thanks Indiana University for
sabbatical support.

  \section{Review of twisted homology and Alexander modules}\label{sec:2}

  In this section we review the basic theory of twisted homology groups.

Let  $X$ be a connected finite CW complex with fundamental group $\pi = 
\pi_1(X)$, let  $R$ be a commutative ring,
and let $M$ be a 
right $R$--module.  Suppose also that $\rho \co \pi  \to 
\mbox{Aut}_R(M)$ is a  homomorphism.  Via $\rho$, $M$ can be viewed as 
a right  $R[\pi]$--module, 
with $m \cdot g = \rho (g^{-1})(m)$ for $m \in M$ and $g \in \pi$.

Let $\tilde{X}$ denote the universal cover of $X$ and let $C_* = 
C_*(\tilde{X})$ denote its cellular chain complex with $R$ 
coefficients.   This  can be viewed as a left $R[\pi]$--module.

Given this, one can form the $R$--chain complex  $M \otimes_{R[\pi]} 
C_*$.  The homology of this complex is denoted $H_*(X, M)$.  Note that 
for  a ring $S$, if $M$ has a left $S$--module  structure which is 
compatible with $\rho$ (meaning $s( m\cdot g) = (sm) \cdot g$)
then $H_*(X,M)$ is naturally a left $S$--module.

Two basic examples the reader should consider are the case that $M = R$ 
with trivial $\pi$ representation, in which case $H_*(X, M) = H_*(X, 
R)$, and the case that $M = R[\pi]$ with standard right $\pi$  action, 
in which case $H_*(X, M) = H_*(\tilde{X}, R)$.
Generalizing   these examples, one has the following basic theorem in homological algebra, Shapiro's Lemma.    See Brown~\cite{br} for a basic reference.

\begin{theorem}[Shapiro's Lemma]  If $M$ is a $\pi$--module  and $\kappa 
\subset \pi$, then $H_*(X_\kappa, M) = H_*(X, M \otimes_R  R[\pi/ 
\kappa])$, where $X_\kappa$ is the cover of $X$ associated to the subgroup $\kappa$ and  $M \otimes_R  R[\pi/\kappa]$ is a $\pi$--module via the 
diagonal action.
\end{theorem}

 \subsection{Examples:\qua  Link exteriors and twisted Alexander modules}  For a link $L \subset S^3$ 
the complement $X_L$ 
has the homotopy type of a CW complex 
 with one  0--cell, $g$ 
1--cells, and $(g-1)$ 2--cells.  This corresponds naturally to a 
presentation of the fundamental group of the complement, 
$\pi = \langle x_1, \ldots , x_g \mathop{|} r_1, \ldots , r_{g-1} \rangle$.

  The chain complex of the universal cover of $X_L$  as an 
$R[\pi]$--module is given by $$0 \to   R[\pi]^{g-1} 
\overset{\partial_2}{\to} R[\pi]^g  \overset{\partial_1}{\to} R[\pi] 
\to 0.$$  The boundary map $\partial_1$ is given by right 
multiplication by the transpose of the vector $[ 1 - x_1, \ldots , 1 - 
x_g]$.  The boundary map $\partial_2$ is given by right  
multiplication 
by a $(g-1) \times g$ matrix called the Fox matrix  with 
entries the Fox derivatives,  $\partial r_i / \partial x_j \in R[\pi]$.

To define a twisted Alexander polynomial of a link, one begins with a representation $\rho\co \pi \to GL_n(R)$.  Fixing any set of meridians to distinct components, $\{m_1, \cdots , m_r\}$, there is also an 
action of $\pi$ on $R[t_1^{\pm1} , \cdots , t_r^{\pm1} ]$, where $m_i$ acts by multiplication by $t_i$. Then $R^n \otimes R[t_1^{\pm1} , \cdots , t_r^{\pm1} ]$ is naturally a $\pi$--module, which we abbreviate $R^n[t_i^{\pm 1}]$.   (In 
the case of greatest interest here, $L$ has one component or $L = K \cup A$ and only a meridian of $K$ is selected, in which case we 
drop the subscript and work with $R[t^{\pm 1}]$.)  

The homology of $X_L$ with coefficients in $R^n[t_i^{\pm 1}]$ is given by a sequence
\begin{equation}\label{c1} 
0  {\to} (R^n[t_i^{\pm 1}])^{g-1} \overset{\beta_2}{\to}  (R^n[t_i^{\pm 1}])^g \overset{\beta_1}{\to}   R^n[t_i^{\pm 1}]  \to 0  
\end{equation} 
where $\beta_2$ and $\beta_1$ are given by $n(g-1) \times ng$ and $ng \times n$ matrices with entries in $R [t_i^{\pm 1}]$.

\begin{definition}The homology groups of this chain complex, the twisted Alexander modules of $L$, are denoted $A_k(L,\rho)(t_i)$.  These are $R[t_i^{\pm 1}]$--modules. The number of variables and the corresponding meridians are identified in context and not notated.
\end{definition}

\section{Torsion invariants, twisted polynomials, and computations}\label{sec:3} 

\subsection{Torsion invariants}  In general, let $S$ be a commutative ring with unity and suppose that  $M$ is a torsion 
$S$--module with presentation $$ S^m \overset{A}{\to} S^n \to M \to 
0,$$ where $A$ is an $m \times n$ matrix with entries in $S$, $m \ge n$.  The ideal 
generated by all $n \times n$ minors of $A$ is an invariant of $M$, denoted $E_0(M)$.  (Note that $A$ is not required to be injective.)
Assuming that $S$ is a Noetherian UFD (for instance, when $S$ is a polynomial ring over $\zz$, $\qq$, or $\zz / p\zz$),   this ideal is contained in a unique minimal
principal ideal called the order of the ideal, with generator  denoted Ord($M$).  If $M$ is not torsion, we set $\mbox{Ord} (M ) = 0$.  There is the following result; see Theorem 3.12
of Hillman~\cite{hi3} for details. 

\begin{theorem}  If $0 \to M \to N \to P \to 0$ is a short exact 
sequence of $S$--torsion modules, then \rm $\mbox{Ord}( M )\mbox{Ord}( P 
) = \mbox{Ord}( N)$.
\end{theorem}

\subsection{Twisted Alexander polynomials} In the case of greatest interest to us, the twisted Alexander modules of a link, we have:

\begin{definition} $\Delta^k_{L,\rho}(t_i) = \mbox{Ord}(A_k(L,\rho)(t_i) ) $.  In the case that $k=1$ we  drop the
 superscript,  $ \Delta_{L,\rho}(t_i) = \Delta^1_{L,\rho}(t_i)$.
\end{definition}

\subsection{Computations}  In general we are not given a 
presentation of $A_k(L,\rho)(t_i) $, but rather the chain complex that
collectively determines these modules.  Furthermore, the matrices are
quite large and the number of minors that arise is larger still.
Computations are simplified by the following observations, basically
Lemma 4.2 of Kirk--Livingston~\cite{kl1}.

\begin{lemma}Suppose that one has a chain complex $\calc$  of abelian groups:
$$ 0 \to A \overset{a \oplus b}{\to} B \oplus C \overset{c + d}{\to} D \to 0,$$with $c$ injective.  Then there is a naturally defined exact  sequence: 
\begin{equation} \label{c2}
0 \to H_1(\mathcal{C}) \to \mbox{coker}(b) \to \mbox{coker}(c) \to H_0(\mathcal{C})  \to 0.
\end{equation}
\end{lemma}

\begin{proof}  We offer here the following simplification of the argument in~\cite{kl1} which was presented only in the case that $R$ is a field.

Consider the  exact sequence of chain complexes $0 \to \calc' \to \calc \to \calc'' \to 0$ with the complexes corresponding to the columns of the following commutative diagram, with all the maps coming from $a, b, c, d$, and the natural inclusions and projections.
$$
   \begin{diagram}\dgARROWLENGTH=1.0em\relax
      \node[2]{0} \arrow{s}{}  \node{0} \arrow{s}{} \node{0} \arrow{s}{} \\
       \node{0} \arrow{e} {}  \node{0} \arrow{e} {} \arrow{s} {}  \node{A} \arrow{e}{ } \arrow{s,r} {}   \node{A} \arrow{e} {}\arrow{s} {}   \node{0}  \\
       \node{0} \arrow{e} {}  \node{B} \arrow{e} {} \arrow{s,r }{}  \node{B \oplus C} \arrow{e} {}\arrow{s} {}   \node{C} \arrow{e} {}\arrow{s} {}   \node{0}  \\
            \node{0} \arrow{e} {}  \node{D} \arrow{e} {} \arrow{s} {}  \node{D} \arrow{e} {}\arrow{s} {}   \node{0} \arrow{e} {}\arrow{s} {}   \node{0}  \\
              \node[2]{0}     \node{0}   \node{0} 
   \end{diagram}
$$
The lemma now follows from the associated long exact sequence.
\end{proof}

We apply this in the case of the sequence given in 
\ref{c1}, rewritten as
$$0   \to   (R[t_i^{\pm 1}])^{n(g-1)}   \overset{a \oplus b}{\to}   (R[t_i^{\pm 1}])^n  \oplus (R[t_i^{\pm 1}])^{n(g-1)}    \overset{c + d}{\to}   R [t_i^{\pm 1}] ^n \to 0,$$  where the choice of splitting of the middle group depends on a choice of an  initial generator of $\pi$.

Applying the previous lemma and recognizing that the orders of the cokernels are given by determinants, we have 

\begin{theorem} \label{wthm} If $A_1(L,\rho)$ and $A_0(L,\rho)$ are torsion, then $$\Delta_{L,\rho}(t_i)/ \Delta^0_{L,\rho}(t_i) = \mbox{det}(b) / \mbox{det}(c).$$
\end{theorem}

In Wada~\cite{wa} this quotient of determinants is defined to be an
invariant of a link, now called the Wada invariant $W$, and is
explicitly shown to be independent of the choice of generator of $\pi$
selected as the initial generator, and more generally, to be
independent of the choice of presentation.  The homological
interpretation given here comes from~\cite{kl1}, and offers an
alternative proof of this independence.

   \section{Relating twisted knot and link polynomials}\label{sec:4}

In our application to periodic knots, one starts with a possibly periodic knot $K$ and hypothesizes an axis $A$ for the action.  One is thus forced to consider the relationship between the twisted invariants of $K$ and those of the hypothesized link $L = K \cup A$.  Specifically, for the knot $K$ there is a representation $\rho\co \pi \to GL_n(R)$ and the   twisted invariants, $A_k(K\!,\rho)(t)$ and 
$\Delta^k_{K\!,\rho}(t)$.  There is the associated link $L$ and the induced representation which we continue to denote by $\rho$.  There is also an induced action of $\pi_1(X_L)$ on $R^n[t^{\pm 1}]$, with the meridian of $A$ acting trivially, and thus there are twisted invariants  $A_k(L,\rho)(t)$ and $\Delta^k_{L,\rho}(t)$.

Before relating these twisted invariants of $K$ and $L$, 
we need a bit of notation.  From the chain complex we have that 
 $$ A_0 (K,\rho)= R^n[t^{\pm 1}]/\langle \{ (I - \rho^*(g))m | g \in \pi, m \in R^n[t^{\pm 1}] 
\}\rangle,$$  
where $\rho^{*}(g) =t^{\epsilon(g)} \rho(g)$ and $\epsilon$ is
the abelianization map $\pi \to  \zz$.  For notational simplicity, we abbreviate this quotient as 
$R^n[t^{\pm 1}] / \mbox{Im}(I - \rho^*)$.  Also, we use the following:

\begin{definition}  $\delta_{L,\rho}(t)  = \mbox{det}(I - \rho(A)t^\lambda) $,
where $\lambda={\rm lk}(K,A)$ denotes the linking number.
\end{definition}

\begin{lemma} $H_1(A, R^n[t^{\pm1}]) = 0$ and $H_0(A, R^n[t^{\pm1}]) $ is presented by $  I- \rho(A)t^\lambda  $ and thus has order $\delta_{L,\rho}(t) $.
\end{lemma}

\begin{proof} Since $A \cong S^1$, it has one 0--cell and one 1--cell.  The action of $\pi_1(A)$ on $R^n[t^{\pm1}]$ is generated by  $  \rho(A)t^\lambda$, and thus the chain complex used to compute its homology is
$$0 \to R^n[t^{\pm1}] \xrightarrow{I - \rho(A)t^\lambda} R^n[t^{\pm1}]
\to 0.$$ As we note at the end of the proof of the next result,
$\lambda \ne 0$, so det$(I - \rho(A)t^\lambda)$ has constant term 1
(set $t=0$), and in particular is nonzero.  Thus the map $ I -
\rho(A)t^\lambda$ is injective and the result follows.
\end{proof}

\begin{theorem}\label{mvthm} In the situation described above:
\begin{enumerate}[\rm(a)]
\item   $A_0(L, \rho) = A_0(K,\rho) = R^n[t^{\pm 1}] / \mbox{{\rm Im}}(I - 
\rho^*)$, and $\Delta^0_{K\!, \rho}(t)  =  \Delta^0_{L, \rho}(t)  $.

\item    If $A_2(K, \rho) = 0$, then there is an exact sequence $$0 \to 
  R^n[t^{\pm 1}] / \langle(I - \rho(A)t^\lambda )(m) | m \in R^n[t^{\pm 1}] \rangle \to  
A_1(L, \rho) \to A_1(K, \rho) \to 0, $$ where $\lambda = \mbox{lk} ( K 
, A)$.

\item    If $A_2(K, \rho) = 0$, then $\Delta_{L, \rho}(t) = 
\delta_{L,\rho}(t)\Delta_{K\!, \rho}(t)$.

\item  If $A_2(K, \rho) = 0$, then $A_2(L, \rho) = 0$.
\end{enumerate}
\end{theorem}

\begin{proof} Statement (a) follows from the definition of the twisted 
homology groups and that the image of $\rho^*$ is the same, whether 
applied to $\pi_1(X_K)$ or $\pi_1(X_L)$.

For part (b), the Mayer--Vietoris sequence associated to the 
decomposition $X_K = X_L \cup N(A)$ yields the following, in which the 
twisted coefficients have been left out of the notation.
\begin{gather*}\hskip-.8in 0 \to H_1(T^2) \to H_1(A) 
\oplus H_1(X_L) \to H_1(X_K) \to \\
\hskip.8in   H_0(T^2)  \to H_0(A) \oplus H_0(X_L) \to H_0(X_K) \to 
0.\end{gather*}
Since the representation is trivial on the meridian of $A$, the map 
$H_0(T^2) \to H_0(A)$ is an isomorphism. This, along with the explicit 
calculation  that  $H_1(T^2) =  R^n[t^{\pm 1}] /$ $\langle(I - \rho(A)t^\lambda  )(m) | m \in R^n[t^{\pm 1}] \rangle $ and that $H_1(A) = 0$,  yields the next 
sequence.
$$ 0 \to R^n[t^{\pm 1}] / \langle I - \rho(A)t^\lambda  \rangle \to H_1(X_L) 
\to H_1(X_K) \to 0. $$
Part (c) follows from (b).

For part (d), from the Mayer--Vietoris sequence 
$$  0 \to H_2(T^2) \to H_2(A) \oplus H_2(X_L) \to H_2(X_K) = 0,$$
we   conclude, since $H_2(A) = 0$ ($A$ is   a 1--complex), $H_2(X_L) = H_2(T^2)$.  From an explicit calculation we have $H_2(T^2) = \mbox{ker} ( I - \rho(A)t^\lambda)$. But clearly the determinant of $I - \rho(A)t^\lambda$ is a nonzero polynomial: as we recall in the next paragraph,   $\gcd(\lambda,q) = 1$, so $\lambda \ne 0$ and the constant term in $\mbox{det}(I  - \rho(A)t^\lambda)$ is 1.

To conclude, we need to explain why $\gcd(\lambda,q) = 1$.  The knot $K$ is the preimage of $\bar{K}$ in the $q$--fold cyclic cover of $S^3 -\bar{A}$.  Since $\bar{K}$ is connected, the number of components of $K$, 1, is given by $\gcd (q, \mbox{lk}(\bar{K},\bar{A}))$.  The cyclic cover can be constructed using the Seifert surface for $\bar{A}$, and $\mbox{lk}(\bar{K},\bar{A})$ is the algebraic intersection of $\bar{K}$ with that surface.  Similarly, $\lambda = \mbox{lk}( {K}, {A})$ is computed as the algebraic intersection of $K$ with a lift of the Seifert surface, which is clearly the same number.
\end{proof}


\section{Application of Shapiro's Lemma to periodic knots}\label{sec:5}

In this section we work with 
a field  $R$, which is 
either $\qq[\zeta_q]$ or $\zz / p\zz$, 
where $\zeta_q$ is a fixed primitive $q$--root of unity and $p$ is prime.

Suppose that $K$ is invariant under a $\zz / q\zz$ action on $S^3$ with 
fixed point set a circle $A$ disjoint from $K$.  We can form the 
quotient link  $\bar{L} = (\bar{K}, \bar{A})$.  Recall that $ {\rho} \co \pi_1(X_K) \to  
GL_n(R)$ and that $R^n[t^{\pm1}] \cong R^n \otimes R[t^{\pm1}]$ is viewed as a   $\pi_1(X_{K})$--module.  In addition, assume 
that there is a representation $\bar{\rho} \co \pi_1(X_{\bar{K}}) \to  
GL_n(R)$ of which our original $\rho$ is the lift.  Shapiro's 
Lemma yields the following.

\begin{theorem} \label{shapapp} $H_1(X_L, R^n[t^{\pm 1}] ) = H_1(X_{\bar{L}}, R^n[t^{\pm 1}] \otimes  
R[\zz / q\zz])$.
\end{theorem}
\noindent Here the module $R^n[t^{\pm 1}] \otimes  R[\zz / q\zz]$ is acted on 
diagonally by $\pi_1(X_{\bar{L}})$, the meridian of $\bar{K}$ acting 
trivially on $R[\zz / q\zz]$ and the meridian of $\bar{A}$ acting by 
multiplication by a fixed generator of $\zz / q\zz$ on $R[\zz / q\zz]$.

\subsection{Application 1:\qua  $R = \qq[\zeta_q]$}

In this case there is a splitting of the $\zz / q\zz$--module:  $R[\zz / q\zz]  = 
\oplus _{i=0}^{q-1}  \qq[\zeta_q^i]$, where  $\qq[\zeta_q^i]$ is  viewed as an $R[\zz / q\zz]$--module, with the generator of 
$\zz / q\zz$ acting on $\qq[\zeta_q^i]$ by multiplication by $\zeta_q^i$.   
Hence, we have the next result.

\begin{theorem}\label{factorthm}    $H_1(X_L, R^n[t^{\pm 1}] )\ = 
\oplus _{i=0}^{q-1} 
H_1(X_{\bar{L}}, R^n[t^{\pm 1}] \otimes  \qq[\zeta_q^i])$.
\end{theorem}

{\bf Remark}\qua Note  that  $R^n[t^{\pm 1}] \otimes \qq[\zeta_q^i]  \cong R^n[t^{\pm 1}]$ with the $ \zz / q\zz$   action  given by multiplication by  $\ \zeta_q^i$.  In particular, $H_1(X_{\bar{L}}, R^n[t^{\pm 1}] 
\otimes \qq[\zeta_q^0]) = H_1(X_{\bar{L}}, R^n[t^{\pm 1}])$ with trivial $\zz / q\zz$ action.

\subsection{Application 2:\qua $R = \zz / p\zz$ and $q = p^r$ for some prime $p$.}

In this case, $R[\zz / q\zz] =  (\zz / p\zz)[g]/\langle 1 - g^q \rangle$, where $(\zz / p\zz)[g]$ is the polynomial ring in the variable $g$.  Since the coefficients are in $\zz / p\zz$, this equals  $(  \zz / p\zz)[g]/\langle (1 - g)^q \rangle$.  In general, write $ ( \zz / p\zz)[g]/\langle (1 - g)^k \rangle = V_k$ and note that $V_1 = \zz / p\zz$. The following sequence of modules is exact:
$$ 0 \to  V_{k -1} \to V_k \to V_1  \to 0,$$ where the map $V_{k-1} \to V_k$ is given by multiplication by $1-g$.
\fullref{shapapp} gives that    $H_1(X_L, R^n[t^{\pm 1}] ) = H_1(X_{\bar{L}}, R^n[t^{\pm 1}] \otimes  
V_q)$.  We now have the following result, a ``splitting theorem'' for the twisted homology, that will provide the basis for an inductive argument computing the twisted Alexander polynomial for $\zz / p\zz$--coefficients in 
\fullref{sectiontwistp}.

\begin{theorem}\label{factorthmp} If $H_2(X_K, R^n[t^{\pm 1}]) = 0$, then for $k>1$ there is the  following exact sequence:
$$\hskip-.7in  0 \to   H_1(X_{\bar{L}}, R^n[t^{\pm 1}] \otimes  
V_{k-1})  \to H_1(X_{\bar{L}}, R^n[t^{\pm 1}] \otimes  
V_k) \to $$
$$ \hskip1in H_1(X_{\bar{L}}, R^n[t^{\pm 1}])  \to H_0(X_K,R^n[t^{\pm1}]) \to    0 .$$ 
\end{theorem}
\begin{proof}  Note first that $   H_1(X_{\bar{L}}, R^n[t^{\pm 1}] \otimes  
V_1) = H_1(X_{\bar{L}}, R^n[t^{\pm 1}])$, since $V_1 = \zz / p\zz$.

The tensor products are taken over $\zz / p\zz$, and $R^n[t^{\pm 1}]$ is free as a $\zz / p\zz$ vector space, so  tensoring with it preserves exactness and the following sequence is exact:
$$ 0 \to  R^n[t^{\pm 1}]  \otimes V_{k -1} \to  R^n[t^{\pm 1}] \otimes V_k \to  R^n[t^{\pm 1}] \otimes V_1  \to 0.$$

{\bf Observation 1}\qua  $H_2(X_{\bar{L}}, R^n[t^{\pm 1}] \otimes V_1) = 0$. 

Since $H_2(X_K, R^n[t^{\pm 1}]) = 0$, by \fullref{mvthm} we have that $H_2(X_L, R^n[t^{\pm 1}]) = 0$.  By
  Shapiro's Lemma, we conclude that
  $H_2(X_{\bar{L}}, R^n[t^{\pm 1}] \otimes ( \zz / p\zz)[g]/ \langle 1 - g^q \rangle) = 0$.  Since the coefficients are in $\zz / p\zz$ this is the same as  $H_2(X_{\bar{L}}, R^n[t^{\pm 1}] \otimes V_q) = 0$.

 Consider now the exact sequence coming from the exact sequence of coefficients: 
  $$ H_3(X_{\bar{L}}, R^n[t^{\pm 1}] \otimes  
V_{1}) \to   H_2(X_{\bar{L}}, R^n[t^{\pm 1}] \otimes  
V_{k-1})  \to H_2(X_{\bar{L}}, R^n[t^{\pm 1}] \otimes  
V_k).$$
Since we are working with homotopy  2--complexes, the first term is 0 and thus for all $k$,  $H_2(X_{\bar{L}}, R^n[t^{\pm 1}] \otimes  
V_{k-1})  \to H_2(X_{\bar{L}}, R^n[t^{\pm 1}] \otimes  
V_k) $ is injective. Starting with the fact that $H_2(X_{\bar{L}},   R^n[t^{\pm 1}] \otimes  
V_{q}) = 0$ and working down, we conclude that $H_2( X_{\bar{L}}, R^n[t^{\pm 1}] \otimes  
V_{1}) = 0$, as desired.

 \medskip
{\bf Observation 2}\qua  $H_0(X_{\bar{L}}, R^n[t^{\pm 1}] \otimes V_{k-1}) \to  H_0(X_{\bar{L}}, R^n[t^{\pm 1}] \otimes V_{k})   $ is the zero map.

 In general, $H_0(X_{\bar{L}}, R^n[t^{\pm 1}] \otimes V_{k})$ is given as the quotient of  
 $ R^n[t^{\pm 1}] \otimes V_{k}$ by the ideal generated by all elements of the form $(I - \bar{\rho})(\alpha)(x \otimes 1)$, where $\alpha \in \pi_1(X_{\bar{L}})$, $x \in R^n[t^{\pm 1}]$, and $\bar{\rho}\co \pi_1(X_{\bar{L}}) \to \mbox{Aut}(R^n[t^{\pm 1}] \otimes V_{k})$, the   action defined earlier on $R^n[t^{\pm1}]$ extended as described to $R^n[t^{\pm 1}] \otimes V_{k}$. 

 A meridian to $\bar{A}$ acts trivially on $R^n[t^{\pm 1}]$ and acts by  multiplication by $g$ on $V_k$.  Thus, in the quotient space, $1-g$ acts trivially.  On the other hand, the map from $V_{k-1}$ to $V_k$ is given by multiplication by $1-g$, and hence the induced map on homology is trivial.
 
 \medskip{\bf Observation 3}\qua There exists a natural isomorphism, $$H_0(X_{\bar{L}},R^n[t^{\pm 1}] \otimes V_k) \to H_0(X_K, R^n[t^{\pm1}]).$$
 
Since $V_k / \langle 1- g \rangle = \zz / p\zz$, the inclusion $H_0(X_{\bar{L}}, R^n[t^{\pm 1}]  ) \to  H_0(X_{\bar{L}}, R^n[t^{\pm 1}] \otimes V_{k})   $ is an isomorphism.  The inclusion $X_{\bar{L}} \to X_{\bar{K}}$ is surjective on $\pi_1$, so induces an isomorphism on $H_0$ with twisted coefficients.  Similarly, the projection $X_K \to X_{\bar{K}}$ is surjective on $\pi_1$, so this too induces an isomorphism on $H_0$ with twisted coefficients.
\end{proof}

\section{Murasugi's condition for twisted polynomials of periodic 
knots with integral or rational representations}\label{sectionproofrat}

Here we prove our main result concerning integral and rational twisted polynomials of periodic knots: 

\medskip
{\bf \fullref{tthm1}}\qua {\sl Suppose that $K$ is of period $q$, and
$A$, $\bar{K}$, and $\bar{A}$ are as above.  Further, assume that
$\bar{\rho} \co \pi_1(S^3 - \bar{K}) \to GL_n(R)$, with $R = \zz$ or
$\qq$, and $\rho$ is the lift of that representation to $S^3 - K$.  If
$\Delta_{K\!, \rho}(t) \ne 0$, then there is a polynomial $F(t,s) \in
R[t^{\pm 1},s^{\pm 1}]$ such that $$\Delta_{K\!, \rho}(t) =
\Delta_{\bar{K}\!, \bar{\rho}}(t) \prod_{i=1}^{q-1} F(t, \zeta_q^i).$$ }

Note that since $\Delta_{K\!,  \rho}(t)  \ne 0$,  $A_2(K,\rho) = 0$  and \fullref{mvthm} applies.

Let $R = \zz \mbox{ or }  \qq$.  The Alexander polynomial of interest is associated to the homology group  $H_1(X_K, R^n[t^{\pm 1}])$.  We can tensor $R^n$ with $\qq[\zeta_q]$ to form the module $R_\zeta^n$.  (Notice that $R_\zeta = \qq[\zeta_q]$, but writing $R_\zeta$ will reduce confusion.)  Denote the representation obtained by tensoring as $\rho_\zeta$.  Observe that $\Delta_{K\!,\rho}  (t)  = \Delta_{K\!,\rho_\zeta}  (t)$, as follows from the computational formula provided by \fullref{wthm}.

By \fullref{mvthm} (part c), $\Delta_{K\!, \rho_\zeta}(t)\delta_{L,\rho}(t) =  \Delta_{L,\rho_\zeta}(t)$. Next,  
we   apply   \fullref{factorthm}:
 $$\Delta_{L,\rho_\zeta}(t) =  \Delta_{\bar{L},\bar{\rho}_\zeta}(t) \prod_{i=1}^{q-1} \mbox{Ord}H_1(X_{\bar{L}},  R^n[t^{\pm 1}]  \otimes \qq[\zeta_q^i]).$$ 

Again by \fullref{mvthm} (part c),   $\Delta_{\bar{L},\bar{\rho}_\zeta}(t) =       \Delta_{\bar{K},\bar{\rho}_\zeta}(t)  \delta_{\bar{L},\bar{\rho}_\zeta}(t)$.  But note that  $\delta_{\bar{L},\bar{\rho}_\zeta}(t) = \delta_{ {L},\rho}(t)$, since $\lambda$ is the same for $L$ and $\bar{L}$ as is $\rho(A)$ and $\bar{\rho_\zeta}(\bar{A})$.

 Then, canceling the $\delta_\lambda$ terms, we have:
$$ \Delta_{K\!, \rho_\zeta}(t)  = \Delta_{\bar{K}\!,\bar{\rho}_\zeta}(t)  \prod_{i=1}^{q-1} \mbox{Ord}H_1(X_{\bar{L}}, R_\zeta^n \otimes \qq[\zeta_q^i]).$$
Thus, it remains to analyze $\mbox{Ord}H_*(X_{\bar{L}}, R_\zeta^n[t^{\pm 1}] \otimes \qq[\zeta_q^i])$.  Suppose that $\pi_1(S^3 - \bar{L})$ has a presentation with $g$ generators and $g-1$ relations, and assume further that the first of the generators, $x_1$, is a meridian to the axis.  There is an action of $\pi_1(S^3 - \bar{L})$ on $R^n[t^{\pm 1}] \otimes \qq[s^{\pm 1}] = R^n[t^{\pm 1}, s^{\pm 1}]$, where the meridians to $\bar{A}$ act  on  $\qq[s]$   by multiplication by  $s$, and  meridians to $\bar{K}$ act  on  $\qq[s]$   by multiplication by  $1$.  The homology groups of $X_{\bar{L}}$ with these coefficients, the two variable twisted 
Alexander modules of $\bar{L}$, are determined by an exact sequence:
\begin{multline} \label{eqn1} 0 \to  (R^n[t^{\pm 1}, s^{\pm 1}])^{g-1} \xrightarrow{A(s,t) \oplus B(s,t) }  (R^n[t^{\pm 1}, s^{\pm 1}])   \oplus  (R^n[t^{\pm 1}, s^{\pm 1}])^{(g-1)} \\ 
 \xrightarrow{C(t,s) + D(t,s)}    R^n[t^{\pm 1}, s^{\pm 1}] \overset{}{\to}  0.\end{multline}
Replacing $s$ with $\zeta_q^i$ in this sequence (that is, tensoring over $\qq[s]$ with $\qq[\zeta_q^i]$ with $s$ acting on $\qq[\zeta_q^i]$ by multiplication by $\zeta_q^i$) gives the sequence used to compute the homology $H_*(X_{\bar{L}}, R_\zeta^n[{t^\pm1}] \otimes \qq[\zeta_q^i]).$ 

The determinant of $C(t,s)$ is simply $(1-s)^n$, and this is nonzero if $s = \zeta_q^i,  1 \le i < q$. Hence, by  \fullref{wthm}, 
$$\hskip-2in \mbox{Ord}H_1(X_{\bar{L}}, R_\zeta^n[t^{\pm 1}] \otimes \qq[\zeta_q^i])=$$
$$\hskip.7in  \mbox{det}(B(t,\zeta_q^i) )\cdot \mbox{Ord} H_0(X_{\bar{L}}, R^n[t^{\pm 1}] \otimes \qq[\zeta_q^i])/\mbox{det}(C(t,\zeta_q^i)).$$  
However, since   $C(t,\zeta_q^i) = (1-\zeta_q^i) I $  and $ 1-\zeta_q^i$ is nonzero, the homology group $H_0(X_{\bar{L}}, R^n[t^{\pm 1}] \otimes \qq[\zeta_q^i]) = 0$.      Letting $G(t,s) = \mbox{det}(B(t,s))/(1-s)^n$, we now have:

\begin{theorem} In the above situation,  
$$ \Delta_{K\!, \rho_\zeta}(t)  = \Delta_{\bar{K}\!,\bar{\rho}_\zeta}(t) \prod_{i=1}^{q-1} G(t,\zeta_q^i).$$
\end{theorem}

It is not clear initially that $G(t,s)$ is a polynomial, because of the factor $(1-s)^n$ in the denominator.  However    we do have the following fact (see also~\cite{wa}).

\begin{lemma}$G(t,s) \in R[t^{\pm 1}, s^{\pm1}]$. 
\end{lemma}
\begin{proof} If we had made a meridian to $\bar{K}$, say $m$, the first 
generator instead of the meridian to $\bar{A}$, we would have had a similar sequence 
to \ref{eqn1}: 
\begin{multline} \label{eqn2} 0 \to  (R^n[t^{\pm 1}, s^{\pm 1}])^{g-1} \xrightarrow{A^*(s,t) \oplus B^*(s,t) }  (R^n[t^{\pm 1}, s^{\pm 1}])   \oplus  (R^n[t^{\pm 1}, s^{\pm 1}])^{(g-1)} \\ 
 \xrightarrow{C^*(t,s) + D^*(t,s)}    R^n[t^{\pm 1}, s^{\pm 1}] \overset{}{\to}  0.\end{multline}
Now we  can define  $G^*(t,s) =  \mbox{det}(B^*((t,s))/\mbox{det}(I-\bar{\rho}(m)t)$. According to \fullref{wthm}, the choice of initial  generator does not affect the quotient, and thus $G(t,s) = G^*(t,s)$.  Written otherwise, $$\mbox{det}(B(t,s)) \mbox{det}(I-\bar{\rho}(m)t) = \mbox{det}(B^*(t,s)) (1-s)^n.$$  Since $\qq[\zeta_q][t^{\pm 1},s^{\pm 1}]$ is a UFD, $\mbox{det}(B(t,s))$ is divisible by $(1-s)^n$ and the result follows.
\end{proof}

\section{Murasugi's condition for twisted polynomials of periodic 
knots with $\zz / p\zz$ representations}\label{sectiontwistp}
 Here we prove our main theorem regarding twisted polynomials with 
$\zz / p\zz$--coeffic\-ients of periodic knots. 

\medskip
{\bf \fullref{tthm2}}\qua {\sl   For a knot $K$ of period $q = p^r$, $p$ prime, and a  representation $\bar{\rho} \co \pi_1(S^3 - \bar{K}) \to GL_n(\zz / p\zz)$ with lift $\rho$,   if  $\Delta_{K\!,\rho}(t) \ne 0$, then 
 $$\Delta_{K\!,\rho}(t) = \Delta_{\bar{K}\!, \bar{\rho}}(t)^q 
 (  \delta_{L ,  {\rho}}(t) /\Delta^0_{K\!,\rho}(t)) ^{q-1}.$$
}

By   \fullref{factorthmp}, we have that $$\hskip-1.5in \mbox{Ord}H_1(X_{\bar{L}}, R^n[t^{\pm 1}] \otimes V_q) \mbox{Ord}H_0(X_K,R^n[t^{\pm 1}]) = $$
$$\hskip1.5in \mbox{Ord}H_1(X_{\bar{L}}, R^n[t^{\pm 1}] \otimes V_{q-1})\mbox{Ord}H_1(X_{\bar{L}}, R^n[t^{\pm 1}]).$$ 
Dividing by the $H_0$ term and proceeding by induction to the case of $V_{q - (q-1)} = \zz / p\zz$, we find that
 $$\hskip-2.5in \mbox{Ord}H_1(X_{\bar{L}}, R^n[t^{\pm 1}] \otimes V_q) = $$
$$\hskip.5in \mbox{Ord}H_1(X_{\bar{L}}, R^n[t^{\pm 1}] )( \mbox{Ord}H_1(X_{\bar{L}}, R^n[t^{\pm 1}])/  \mbox{Ord}H_0(X_K,R^n[t^{\pm 1}]))^{q-1} .$$ 
Applying Shapiro's Lemma, \fullref{shapapp}, and switching notation,  this gives that
$$\Delta_{L,\rho}(t) = \Delta_{\bar{L},\bar{\rho}}(t)^{q}/\Delta^0_{K\!,\rho}(t)^{q-1}.$$
Since $ \Delta_{K\!,  \rho}(t)  \ne 0$ we have $ A_2(K,\rho) =  0$.  We then have, by \fullref{mvthm}, that $\Delta_{K\!,\rho}(t) \delta_{L,\rho} (t)= \Delta_{L,\rho}(t)$, and similarly for $\bar{K}$ and $\bar{L}$,  with $\delta$ defined as before.
Thus,  
$$\Delta_{K\!,\rho}(t) \delta_{L,\rho} (t) =( \Delta_{\bar{K}\!,\bar \rho}(t) 
\delta_{\bar{L},\bar \rho} (t))^{q}/\Delta^0_{K\!,\rho}(t)^{q-1}.$$
As we have observed earlier, $\delta_L = \delta_{\bar{L}}$, so further simplification completes the proof of the theorem.


    \section{Examples}\label{10_162}
    
   We use the classical numbering system for prime knots of 10 or fewer crossings, eliminating  the duplicate $10_{161} = 10_{162}$ that appeared in early tables.  There are 165 prime knots of 10  crossings.  For prime knots of 11 and 12 crossings we use the numbering system based on Dowker--Thistlethwaite notation~\cite{dt}.  This number system is used by the program Knotscape, developed by Hoste and Thistlethwaite,~\cite{ht}.  An up-to-date table of knots through 12 crossings is available on-line at~\cite{lc}.
    
 Consider the knot $10_{162}$, having irreducible 
Alexander polynomial $$\Delta_{10_{162}}(t) = 3t^4 - 9t^3 +11t^2 - 9 t 
+3. $$
  We first observe that Murasugi's conditions do not obstruct this knot 
from having period 3. 
One simply checks that 
letting   $\Delta_{\bar{K}}(t) = \Delta_{10_{162}}(t)$ and $\lambda = 1$ 
yields a 
solution  to the mod 3 condition of \fullref{murthm2}.

Now, for the $\zz[\zeta_3]$ condition in \fullref{murthm1}, one 
lets $\alpha = 1+s+s^2$ and 
$$F(t,s) = \alpha t^4 - 3\alpha t^3 + (4\alpha -1)t^2 - 3\alpha t +  \alpha.$$
   We want to see a bit more, 
that this is the only possible solution to 
the Murasugi conditions. Suppose that $10_{162}$ did have period 3.  
As $\Delta_{10_{162}}(t) \equiv 1\ {\rm mod\ }3,$ 
using the mod $3$ condition
we see that the linking number must be 
$\lambda = 1$.  Now using the $\zz[\zeta_3]$ criterion
we have that $ \Delta_{10_{162}}(t)$ would factor as $\Delta_{\bar{K}}(t) F(t) 
F(t)^\sigma$, where $F(t) \in \zz[\zeta_3][t, t^{-1}]$ and $\sigma$ is 
a nontrivial Galois automorphism in the Galois group of  $\qq(\zeta_3)$ over $\qq$.  However, $\Delta_{10_{162}}(t)$ is 
irreducible in $\zz[\zeta_3][t, t^{-1}]$, and hence the only factor 
must be $\Delta_{\bar{K}}(t) =  \Delta_{10_{162}}(t)$.

   We now want to observe that the twisted Murasugi conditions do 
obstruct $10_{162}$ from having period 3.  The previous argument shows 
that if it does have period 3, the quotient knot $\bar{K}$ has the same 
polynomial.  Thus,  the 2--fold branched cover of the quotient knot has 
homology of order $\Delta_{10_{162}}(-1) = 35$.  Since this is 
divisible by $5$, but not $5^2$, it follows that $\bar{K}$ has a surjective  representation to the fifth 
dihedral group, $D_{5}$, that is unique up to composition with automorphisms of $D_{5}$.  (Recall that representations to $D_5$ correspond to representations of the first homology of the 2--fold branched cover to $\zz / 5\zz$.)

   The group $D_{5}$ has a representation to $GL_4(\zz)$.    It is 
defined via the identification of $\zz^4$ with 
$\zz[\zeta_5]$, as free abelian groups: 
the  order two automorphism is given by complex 
conjugation and the order 5 automorphism is given by multiplication by 
$\zeta_5$.  Similar reasoning applies to the knot $10_{162}$, and so $10_{162}$ has a nontrivial representation $\rho$ to $GL_4(\zz)$ which lifts a representation $\bar{\rho }$ of $\bar{K}$ with image this dihedral subgroup.   

   With respect to this representation, a computation gives that
    $$ \Delta_{10_{162}, \rho},(t) = (11 t^8-21 t^6-39 t^4-21 t^2+11) 
(t-1)^3 (t+1)^3.$$
   Each factor in the above factorization is irreducible in 
$\zz[\zeta_3][t, t^{-1}]$.  Hence, the $\zz[\zeta_3]$ criterion in
\fullref{tthm1} implies 
that $\Delta_{\bar{K}\!, \bar{\rho}}(t) =  \Delta_{10_{162} , \rho}(t)$.  
Given this, it becomes evident that the mod 3 condition 
in \fullref{tthm2} 
cannot hold.  
Thus, the knot does not have period 3.

This knot was already known not to have period 3: see~\cite{ksa} or~\cite{na}.
It should, however, be noted that by a result of~\cite{dl1}, any Alexander 
polynomial $\Delta$ satisfying $\Delta \equiv \pm t^m$ mod $n$ for 
some $m$ is the Alexander polynomial of some knot of period $n$.

\section{Alexander polynomial 1 examples}
\label{poly1}
 As a second and   more subtle application, we consider the two prime 
knots of 11 crossings having trivial Alexander polynomial:  
$11_{n34}$ and $11_{n42}$.  In both cases we show that the twisted 
Alexander polynomial obstructs periodicity, where clearly the classical 
polynomial cannot.  The  bounds on the periods  for these knots that we 
present here were first found by Naik~\cite{na}, using  geometric work 
of Edmonds~\cite{ed}, based on the theory of minimal surfaces in 
3--manifolds.  
 Alternatively, SnapPea also 
provides the periods
 of these knots.  We present these examples here both to illustrate our 
purely algebraic approach, and to indicate the independence of the 
results from  deep differential geometric theory.

\begin{figure}[ht]
\centerline{  \fig{.35}{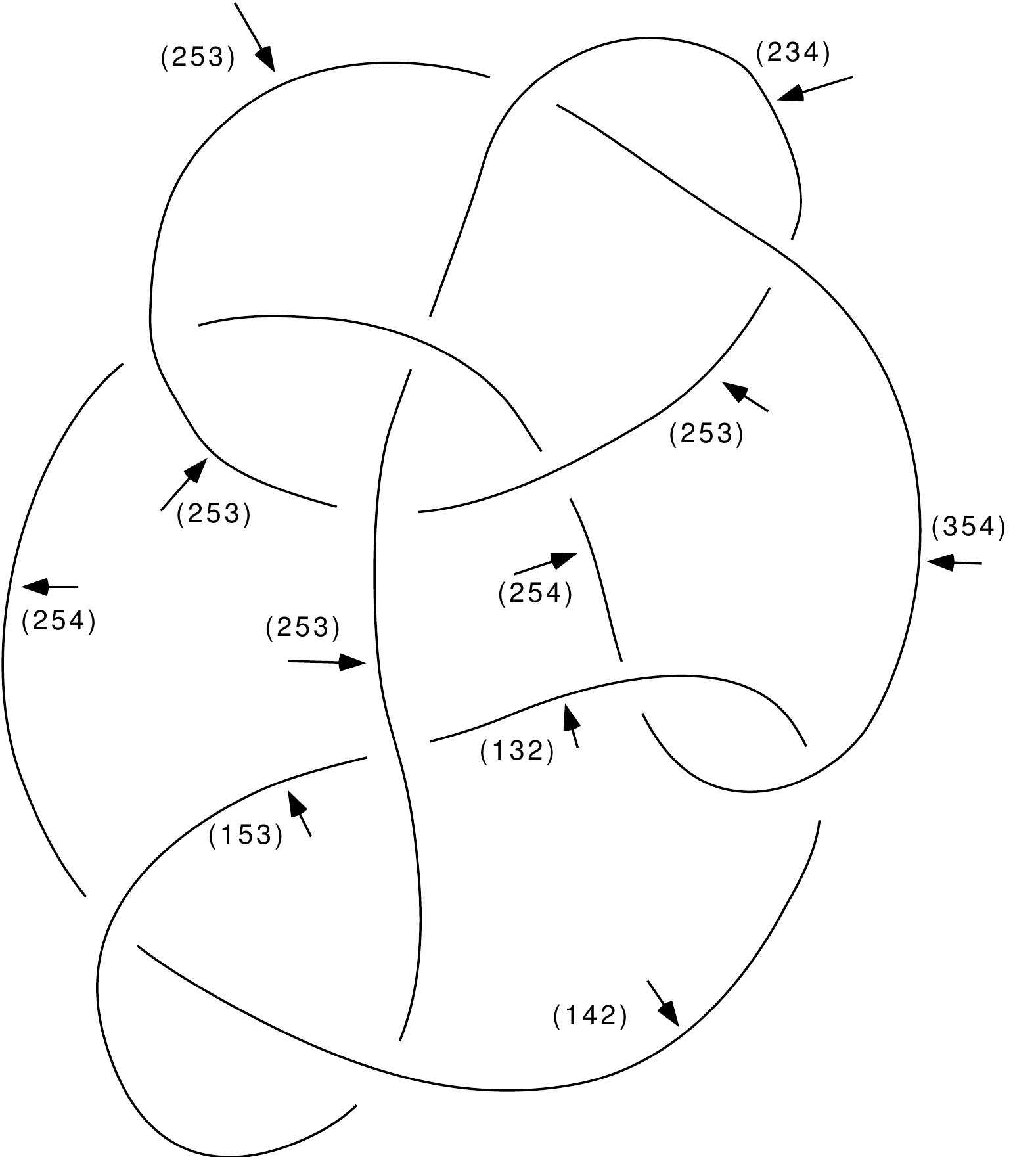}  \hskip.1in  \fig{.35}{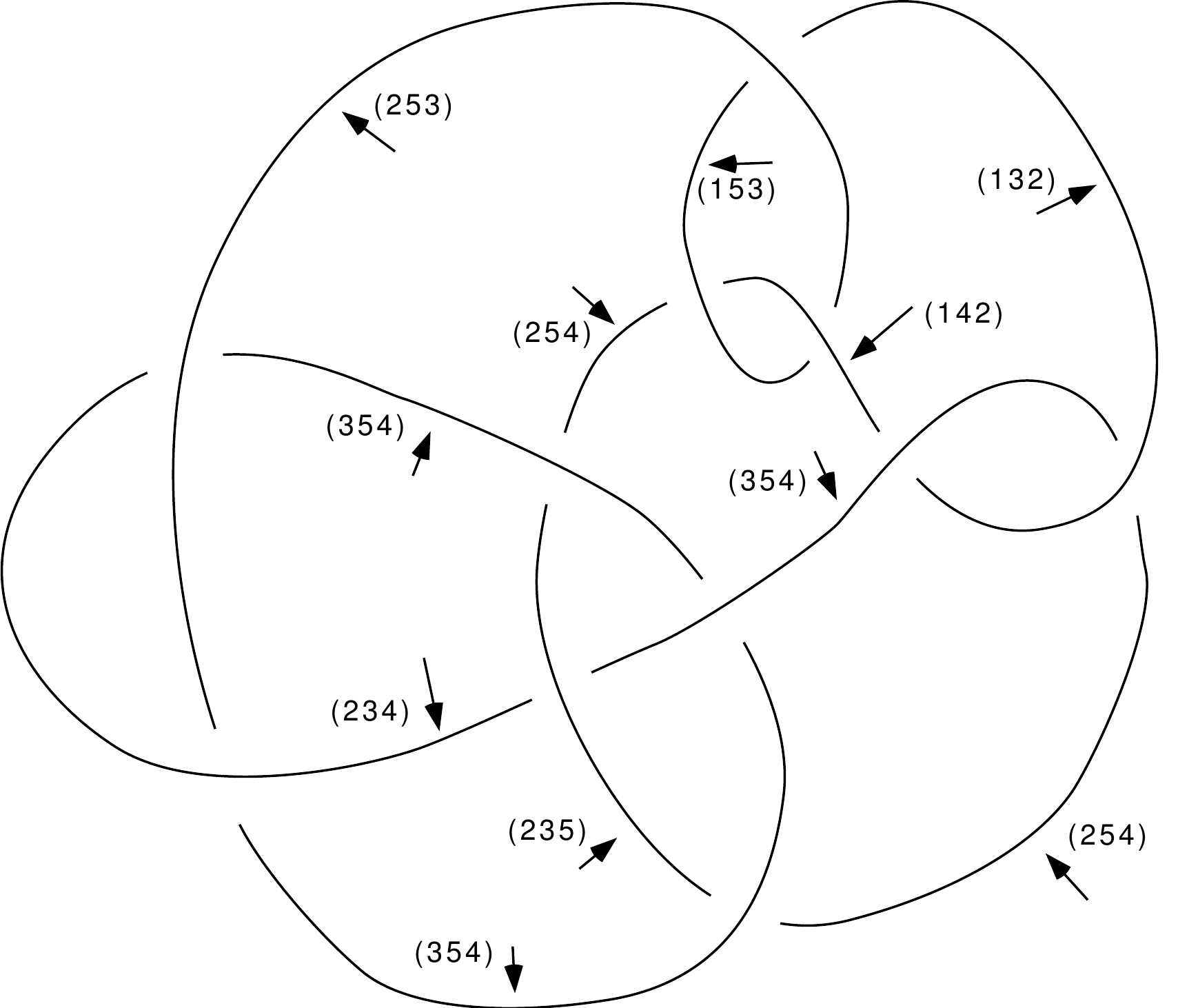}  }
 \caption{$11_{n34}$  and $11_{n42}$} \label{fig10_n3442}
\end{figure}

 Both of these knots have representations to the alternating group, $A_5$.  
In \fullref{fig10_n3442} we have illustrated both knots, along with one 
choice of representations.
 
 Suppose that one of these knots, we call it $K$, is periodic of 
period $q$ with axis $A$.  Selecting a base point for the fundamental group 
on $A$, we have the $\zz / q\zz$  action inducing an action on the 
representations of $\pi_1(X_K)$ onto $A_5$.  However, a direct  calculation 
shows that there is a unique surjective  representation up to an  
automorphism of $A_5$.  (This   calculation was done with 
Knotscape~\cite{ht}.)  Every automorphism of $A_5$ is given 
by conjugation by an element of $S_5$, and in particular, 
the only possible orders of automorphisms are 2, 3, 4, 5, and 6.  
 
 Let $p$ be a prime other than 2, 3, or 5.  It follows that if $K$ has 
period $p$, then the action fixes the representation illustrated in the 
figure, and thus the representation is the lift of a representation 
to $A_5$ of the quotient knot.  
  
\subsection{Example:\qua $11_{n34}$}
To apply Theorems~\ref{tthm1}~and~\ref{tthm2} we use the standard representation of $A_5$ acting on $\zz^5$, and then reduce to either rational or $\zz / p\zz$ coefficients.  We begin with the knot $11_{n34}$, which we abbreviate for now as $K$.  Computations based on \fullref{wthm} and using Maple yield  
$$\Delta_{K\!,\rho}(t) = (t^2 -t +1)(t^2+t+1)(5t^4 +t^3 +8t^2 + t +5)(t+1)^2 (t-1)^4,$$ a polynomial of degree 14.

The exact determination of $\delta_{L ,\rho}$ is not possible, but from the classical Murasugi congruence we have that lk($K,A) = \lambda =1$, and since the representation takes values in $SL_5(Z)$, it follows that $\delta_{L ,\rho}(t)$ is a monic polynomial of degree exactly 5.

Finally, we have by a direct computation that $\Delta^0_{K\!,\rho} = (1-t)$.

Consider  the $\zz / p\zz$ criterion of \fullref{tthm2}.  Since we are 
assuming that $p$ is not 5, $\Delta_{K\!, \rho}$ is degree 14, and 
so we have  $4(p-1) + kp = 14$, where $k$ is the degree of the $\zz / p\zz$ 
polynomial $\Delta_{\bar{K}\!, \bar{\rho}}(t)$.  This gives that $(k+4)p = 
18$, an  impossibility, given that $p$ is a prime greater than 5.

 \subsection{Example:\qua $11_{n42}$}
 An identical calculation applies for $K_2 = 11_{n42}$.  The only 
change throughout is that 
$\Delta_{K_2,\rho}(t) = (t-1)^4(t^2+t+1)(5t^6+5t^5-5t^4-9t^3-5t^2+5t+5)$, 
a polynomial of degree 12.  This now gives that $4(p-1) + kp = 12$, 
or $(k+4)p = 16$.  Again, this is not possible for a prime greater than 5.

\section{Additional conditions on twisted polynomials of periodic knots}
\label{sectionpres}

In the Sections \ref{sectionproofrat} and 
\ref{sectiontwistp} we obtained conditions on the twisted polynomials of 
a periodic knot when the representation of the fundamental group of the
knot complement is a lift of a representation for the
quotient knot. In this section we consider representations
which do not factor through the fundamental group of the quotient.

\begin{theorem}\label{thmpres}
Suppose that the knot $K$ has period $q$. 
Then there is a $\zz _q$ action on the 
 representation space\ ${\rm Hom}(\pi_1(S^3-K),GL_n(R))$
and the 
twisted Alexander polynomials are preserved under this action.
\end{theorem}
\begin{proof} 
Let $f$ be the order $q$
homeomorphism of $S^3$ which leaves $K$ invariant and fixes a 
disjoint circle $A$. Choose as a base point for the knot 
group a point in $A$ and
denote by $f_*$ the automorphism of
the knot group induced by the restriction  of $f$ to $S^3-K$.
Let ${\rho} \co \pi_1(S^3 - {K}) \to  GL_n(R)$ be a representation.
The $\zz _q$ action is given by
$f^*\rho= \rho \circ f_*$. 
Clearly, $\Delta _{K\!, \rho } =
\Delta _{K , f^*\rho } .$
\end{proof}

The usefulness of this observation is illustrated in the example 
below where we 
use homological conditions given by the untwisted 
Alexander polynomial to first conclude that the action on the 
representation space has to be nontrivial and then show that 
the twisted polynomials are not preserved under a nontrivial 
$\zz / q\zz$ action, thereby ruling out period $q$.

\medskip
{\bf Example}\qua Consider the knot $K=12_{n847}$. The Alexander polynomial
for $K$ 
is  $\Delta(t)=1-7t+18t^2-23t^3+18t^4-7t^5+t^6=
(t^2-t+1)(t^2-3t+1)^2$. Over 
$\zz [\zeta _3]$ it factors as 
$(t+\zeta _3)(t+\zeta _3^2)(t^2-3t+1)^2$ and it is
congruent to $(1+t+t^2+t^3)^2$ mod 3. Applying the Murasugi conditions 
in Theorems~\ref{murthm1} and~\ref{murthm2} we see that, 
if this knot has period 3, then the polynomial of the quotient knot
$\bar K$ has to be 
$\bar \Delta (t) =1$.

The homology of the 2--fold  cover $M_2$ of $S^3$ branched over $K$
is $ \zz _5 \times \zz _5 \times \zz _3.$
The group $\zz _5 \times \zz _5$ 
has 24 nontrivial representations to $\zz _5$.
These representations correspond to 6 nontrivial conjugacy classes of  
representations of $\pi_1(S^3 - {K})$ to $D_5$. 
As in \fullref{10_162} we view these as  representations to $GL_4(\zz)$.
 
Suppose that $K$ has period 3. Again choosing a base point in the 
fixed set of the order 
$3$ homeomorphism $f$ of $S^3$, we have a unique lift to
the cover $M_2$ which fixes a chosen lift of the base point.
This lift induces a $\zz _3$ action on the 5--primary
subgroup $\zz _5 \times \zz _5$ of $H_1(M_2)$.
A transfer argument gives us that the fixed point set of this 
action is isomorphic to a subgroup of $H_1$ of the quotient $\bar M_2$,
which is trivial since $\bar \Delta (t) =1$.
(See Proposition 2.5 in \cite{na} and Proposition 4 in \cite{na2}.)
It follows that the action is nontrivial.

Therefore the $\zz _3$ action on the corresponding
$GL_4(\zz)$ representations of the knot group 
is nontrivial and it permutes the 6 conjugacy classes. 
Therefore these  representations should give at most 
two distinct twisted polynomials.
However, explicit  computations have resulted in 
4 distinct twisted polynomials thus ruling out period 3 for this
knot.
 
\section{Freely periodic knots}\label{sec:11}

Historically the main focus of the study of periodicity of knots has concerned periodic actions for which there is an axis, the ``periodic knots'' studied so far in this paper.  It is for these knots that the periodic nature can   most readily be illustrated in a knot diagram.  There is a second class of periodic knots, {\it freely periodic knots}, that are invariant under a free action on the 3--sphere.  These are of independent interest  and also arise in the study of knots in lens spaces.

The main classical result concerning freely periodic knots is Hartley's \cite{ha}:

\begin{theorem}\label{hartleythm} If $K$ is a freely periodic knot of period $q$, then $\Delta_K(t^q) = \prod_{i=0}^{q-1} \bar{\Delta}(\zeta_q^i t)$ for some Alexander polynomial $\bar{\Delta}$.
\end{theorem}

In this section we will   prove a generalization of Hartley's theorem to the case of twisted Alexander polynomials  and illustrate the strength of
the generalization by ruling out free periodicity, of all possible orders, for a 10 crossing knot that Hartley was unable to resolve.  Since Hartley's theorem and its generalization call on factorizations in a cyclotomic ring, as in~\cite{ha} we call on some nontrivial number theory.  The number theoretic techniques we use, however, represent a great simplification over those originally used by Hartley.

 In the case of free periodicity, since  there is no axis,  the quotient of $(S^3, K)$, $(M,\bar{K})$, is a knot in a homology lens space.  
To simplify our considerations  we will work with 
$R=\qq[\zeta_q].$

Observe first that $\pi_1(M) = \zz / q\zz$.  We next want to see that  $H_1(M - N(\bar{K})) \cong \zz$.   Since $M- N(\bar{K})$ has torus boundary, its  homology has rank at least one.  (For a 3--manifold $X$ with boundary $\partial = T^2$, the long exact sequence of the pair $(X, \partial)$ yields $\chi(\partial) - \chi(X) + \chi(X,\partial) = 0$.  But $\chi(\partial) = 0$ and $\chi(X,\partial) = -\chi(X)$ by Poincar\'e duality.  Thus, $\chi(X) = 0$.  Since $H_3(X) = 0$ for a bounded 3--manifold, and $H_0 = \zz$, it must be that the first Betti number is positive.)  Attaching $N(\bar{K})$ onto $M- N(\bar{K})$ adds one relation on homology and yields $\zz / q\zz$,   so $H_1(M - N(\bar{K}))$ must be isomorphic to either $\zz$ or $\zz \oplus \zz / q\zz$.   To see that  $H_1(M - N(\bar{K}))$ is torsion free, we recall first that  $\bar{K}$ lifts to a connected curve in the $q$--fold cover of $M$.  It follows  that the inclusion $H_1(\partial N(\bar{K})) \to H_1(M - N(\bar{K}))$ is surjective.  Thus $H_1(M - N(\bar{K}), \partial N(\bar{K})) = 0$.  By duality this implies that $H^2(M - N(\bar{K}) )= 0$.  Finally the universal coefficient implies that $H_1(M - N(\bar{K}))$ is torsion free.  In summary, we have the following lemma.

\begin{lemma} $H_1(M - \bar{K}) \cong \zz$ and the covering $S^3 - K \to M - \bar{K}$ is a $q$--fold cyclic cover corresponding to a surjection of $\zz$ onto $\zz / q\zz$.  The map $H_1(S^3 - K) \to  H_1(M - \bar{K}) $ is the homomorphism from $\zz$ to $\zz$ given by multiplication by $q$.
\end{lemma}

As before, we now write $X_K$ and $X_{\bar{K}}$ for the knot complements.

To apply our techniques to twisted polynomials of $K$, we assume the existence of a representation $\bar{\rho}:  \pi_1(X_{\bar{K}} )\to GL_n(R)$.  The lift of the representation to $X_K$ will be denoted $\rho$. 

There is the natural representation $\bar{\epsilon}$ of $\pi_1(X_{\bar{K}}) \to \zz$ and thus a twisted polynomial $\Delta_{\bar{K}\!,\bar{\rho}}$.  This lifts to give a representation $\epsilon$ of $\pi_1(X_K) \to \zz$, but note now that this is $q$ times the usual representation.  If we denote the twisted polynomial of $K$ associated to $\rho$ and $\epsilon$ by $\Delta_{K\!, \rho, \epsilon}$ the next result follows immediately.

\begin{lemma}  $\Delta_{K\!, \rho, \epsilon}(t) =  \Delta_{K\!, \rho}(t^q)$.  
\end{lemma} 

We can now apply Shapiro's lemma to this situation, as was done in the
proof of \fullref{tthm1} in \fullref{sectionproofrat}, to find
$$\Delta_{K\!,\rho}(t^q) =     \prod_{i=0}^{q-1} \mbox{Ord}H_1(X_{\bar{K}},   R^n[t^{\pm 1}]  \otimes \qq[\zeta_q^i]).$$
Here the action of $\pi_1(X_{\bar{K}})$ on $\qq[\zeta_q^i]$ is via
$\epsilon$.  Because of this, $$R^n[t^{\pm 1}] \otimes \qq[\zeta_q^i
]\cong R^n \otimes R[t^{\pm 1}] \otimes \qq[\zeta_q^i ] \cong R^n
\otimes R[t^{\pm 1}] .$$ (Recall that $R = \qq[\zeta_q]$.) In the last
of these, the isomorphism is as abelian groups.  The module structure
has now changed in that the action of $\pi_1(X_{\bar{K}})$ on
$R[t^{\pm 1}]$ is such that if an element in $\pi_1$ generates
$H_1(X_{\bar{K}})$, it acts by multiplication by $\zeta_q^i t$, rather
than simply by $t$.

\fullref{freethm} from the introduction can now  be restated:   

\medskip
{\bf \fullref{freethm}}\qua {\sl In the notation above:
$$\Delta_{K\!,\rho}(t^q) =   \prod_{i=0}^{q-1} \Delta_{\bar{K}\!,\bar{\rho}}(\zeta_q^i t).$$}

In the case of trivial one dimensional representations, this gives a new proof of Hartley's result, \fullref{hartleythm}.

\medskip{\bf Examples}\qua  
Consider the knot $10_{62}$, a knot identified by Hartley 
as having possible free periods that cannot be ruled out by 
\fullref{hartleythm}
nor by the criterion in terms of  homology of the 2--fold branched cover
obtained in Theorem 2.2 of~\cite{ha}.

The Alexander polynomial for this knot is 
$\Delta = t^8-3t^7+6t^6-8t^5+9t^4-8t^3+6t^2-3t+1$, the product of 
cyclotomic polynomials, $\Phi_{10}(t) (\Phi_6(t))^2$. As Hartley notes 
in~ \cite{ha}, both the factoring condition on Alexander polynomials of 
freely periodic knots and the structure of $H_1(M_2)$ 
leave open the   possibility of  
a free period $q$, for any $q$ which is relatively prime to 30, 
and for these one would have $\Delta_{\bar K}=\Delta$.

Observe that $\Delta_{\bar{K}}(-1) = 45$.  Thus, the quotient knot complement 
would have a $\zz _5$ in the homology of the
2--fold cyclic branched cover giving a nontrivial representation of the knot group onto the dihedral group $D_5$ which would lift to a nontrivial representation of the group of $10_{62}$ onto $D_5$. 
A computer calculation of the associated  twisted polynomial yields a degree 28 polynomial with a degree 8 irreducible factor
$f(t)=t^8-3t^6+3t^4-3t^2+1$ with multiplicity 2. 

If $10_{62}$ has free period $q$, then applying \fullref{freethm} we would have
$$f(t^q) = g(t) \prod_{i=1}^{q-1} g(\zeta_q^it),$$ for some rational polynomial $g$.  (Here we have used that $\qq[\zeta_q][t]$ is a UFD to conclude that the 
factorization given in \fullref{freethm} determines a similar 
factorization for each factor which is irreducible over $\qq$.) 
We want to show  that such a factorization cannot exist.  
As we are showing that there are no free periods, it is sufficient to 
rule out periods of prime order, so we assume henceforth that $q$ is prime.  
For $q \le 11$ that $f(t^q)$ has no such factorization can be checked 
via computer, so we assume that $q \ge 13$.

Initially we know only that   $g$ is a rational polynomial and that the 
equality is   true only up to a rational multiple; we now want to 
observe that $g$ can be assumed to be integral and that the equality is exact.  
If we replace $g$ with its monic associate  we see the equality becomes 
necessarily exact.  Note that evaluating the equation above at $t=0$ 
yields that $g$ has constant term 1.

We now want to show that $g$ has integer coefficients.  
Each root of $g$ is the root of some irreducible factor of $g$, which itself must be a factor of $f$.  But the monic irreducible factors of a monic integral polynomial are integral. Thus the product of the monic irreducible factors of $g$ is also integral, and so $g$ is as well.  

There is a homomorphism $\zz[\zeta_q] \to \zz / q\zz$ sending $\zeta_q $ to 1.  Applying this to the polynomial factorization and recalling that $f(t^q) \equiv (f(t))^q \mod q$, we have
$$f(t)^q \equiv (g(t))^q \mod q.$$
From this it follows that 
$$f(t) = g(t) + qh(t)$$
for some integral polynomial $h$.  Since $f$ and $g$ are monic of degree 8, we know that $h$ is of degree 7 or less.  We also know that $1= f(0) = g(0) + qh(0)$ and $g(0) = 1$, so $h(0) = 0$.

 Next observe that $g(1)$ divides $f(1) = -1$, so   $g(1) = \pm 1$.  Since $f(1) = g(1) +qh(1)$, we have that $1 - \pm 1 = qh(1)$, so clearly $h(1) = 0$ and $g(1) = -1$.
 
 Similarly,  $g(-1)$ divides $f(-1) = -1$, so   $g(-1) = \pm 1$.  But $f(-1) = g(-1) +qh(-1)$, so again $h(-1) = 0$ and $g(-1) = -1$.
 
 Next we move to the ring of Gaussian integers, $\zz[i]$,  and observe that $g(i)$ divides $f(i) = 11$.  It is easily checked that $11$ is prime in $\zz[i]$, so we have that $g(i) = \pm 1, \pm i, \pm 11$, or $\pm 11 i$. In none of these cases would $f(i) - g(i)$ be divisible by $q$ in $\zz[i]$ (we are assuming that $q \ge 13$).  Thus $h(i) = 0$.
 
 By conjugation we have $h(-i) = 0$. 
 
Now we move to the ring $\zz[\zeta_3]$.  A computation shows that $f(\zeta_3) = 5 \zeta_3$.  The only units in $\zz[\zeta_3]$ are of the form $ \pm \zeta_3^i$, and we will let $u$ denote a unit of this form.  Then $g(\zeta_3) = u$ or $g(\zeta_3) = 5u$.  One now checks that for none of the possibilities is $f(\zeta_3) - g(\zeta_3)$ divisible by $q$ (mostly simply, none of the norms of each is divisible by the square of a  prime greater than $11$).  Hence, $h(\zeta_3) = 0$.  Conjugating we have also $h(\zeta_3^2) = 0$.

Finally we work in $\zz[\zeta_6] = \zz[\zeta_3]$  ($\zeta_6 = -\zeta_3$).  Now we find that $f(\zeta_6) = -5\zeta_6 = 5\zeta_3$.  As in the previous paragraph it follows that $h(\zeta_6) = h(\zeta_6^5) = 0$.

At this point we have found nine zeros for the degree 7 polynomial $h$, and thus $h$ is identically 0.  We have now that $g(t) = f(t)$, so that $f(t)$ is a factor of $f(t^q)$.  But $f$ has eight distinct roots and  the one with largest norm, say $\alpha$, is  real, with value approximately 1.47.  This leads to a contradiction: since $f(t) $ divides $f(t^q)$, we have that $f(\alpha^q)=0$.  But $\alpha^q$ has norm greater than any roots of $f$.  

Thus, we have proved that $10_{62}$ does not have any free periods.

Note that all free periods for knots up to 10--crossings not settled by
Hartley were determined by the combined works of Boileau and
Zimmermann in \cite{boz} and Sakuma in \cite{sa2}.  In particular, it
was shown in \cite{boz} that $10_{62}$ is not freely periodic.  The
proofs in \cite{boz} and \cite{sa2} strongly rely on the hyperbolic
structure on the complement unlike our proof above which is based only
on homological algebra.

Also note that Chbili in~\cite{chb} obtains criteria for free periodicity in
terms of the HOMFLYPT polynomial. These criteria are then applied to 
examples in~\cite{ha} for which 
free periods could not be eliminated using \fullref{hartleythm} alone.
Free periods for $K=10_{62}$ are not ruled out in \cite{chb}.

\bibliographystyle{gtart} \bibliography{link}

\end{document}